
\def\LaTeX{%
  \let\Begin\begin
  \let\End\end
  \let\salta\relax
  \let\finqui\relax
  \let\futuro\relax}

\def\UK{\def\our{our}\let\sz s}
\def\USA{\def\our{or}\let\sz z}



\LaTeX

\USA


\salta

\documentclass[twoside,12pt]{article}
\setlength{\textheight}{24cm}
\setlength{\textwidth}{16cm}
\setlength{\oddsidemargin}{2mm}
\setlength{\evensidemargin}{2mm}
\setlength{\topmargin}{-15mm}
\parskip2mm


\usepackage{amsmath}
\usepackage{amsthm}
\usepackage{amssymb}
\usepackage[mathcal]{euscript}

\usepackage[usenames,dvipsnames]{color}
%
%


\def\pier{\color{blue}}
\let\pier\relax
\def\juerg{\color{red}}
\let\juerg\relax


\bibliographystyle{plain}


%

\finqui

\def\Beq{\Begin{equation}}
\def\Eeq{\End{equation}}
\def\Bsist{\Begin{eqnarray}}
\def\Esist{\End{eqnarray}}

\def\Bthm{\Begin{theorem}}
\def\Ethm{\End{theorem}}
\def\Blem{\Begin{lemma}}
\def\Elem{\End{lemma}}

\def\Bcor{\Begin{corollary}}
\def\Ecor{\End{corollary}}
\def\Brem{\Begin{remark}\rm}
\def\Erem{\End{remark}}

\def\Bdim{\Begin{proof}}
\def\Edim{\End{proof}}
\let\non\nonumber




\def\step #1 \par{\medskip\noindent{\bf #1.}\quad}


\def\Lip{Lip\-schitz}
\def\holder{H\"older}
\def\aand{\quad\hbox{and}\quad}

\def\lhs{left-hand side}
\def\rhs{right-hand side}
\def\sfw{straightforward}


\def\generaliz{generali\sz}

\def\organiz{organi\sz}

\def\regulariz{regulari\sz}


\def\multibold #1{\def\arg{#1}%
  \ifx\arg\pto \let\next\relax
  \else
  \def\next{\expandafter
    \def\csname #1#1#1\endcsname{{\bf #1}}%
    \multibold}%
  \fi \next}

\def\pto{.}

\def\multical #1{\def\arg{#1}%
  \ifx\arg\pto \let\next\relax
  \else
  \def\next{\expandafter
    \def\csname cal#1\endcsname{{\cal #1}}%
    \multical}%
  \fi \next}


\def\multimathop #1 {\def\arg{#1}%
  \ifx\arg\pto \let\next\relax
  \else
  \def\next{\expandafter
    \def\csname #1\endcsname{\mathop{\rm #1}\nolimits}%
    \multimathop}%
  \fi \next}

\multibold
qwertyuiopasdfghjklzxcvbnmQWERTYUIOPASDFGHJKLZXCVBNM.

\multical
QWERTYUIOPASDFGHJKLZXCVBNM.

\multimathop
dist div dom meas sign supp .


\def\accorpa #1#2{\eqref{#1}--\eqref{#2}}
\def\Accorpa #1#2 #3 {\gdef #1{\eqref{#2}--\eqref{#3}}%
  \wlog{}\wlog{\string #1 -> #2 - #3}\wlog{}}


\def\graffe #1{\mathopen\{#1\mathclose\}}

\def\<#1>{\mathopen\langle #1\mathclose\rangle}
\def\norma #1{\mathopen \| #1\mathclose \|}

\def\iot {\int_0^t}
\def\ioT {\int_0^T}
\def\iO{\int_\Omega}
\def\intQt{\int_{Q_t}}
\def\intQ{\int_Q}

\def\dt{\partial_t}
\def\dn{\partial_\nu}

\def\cpto{\,\cdot\,}

\def\checkmmode #1{\relax\ifmmode\hbox{#1}\else{#1}\fi}
\def\aeO{\checkmmode{a.e.\ in~$\Omega$}}
\def\aeQ{\checkmmode{a.e.\ in~$Q$}}

\def\aat{\checkmmode{for a.a.~$t\in(0,T)$}}


\def\erre{{\mathbb{R}}}
\def\enne{{\mathbb{N}}}




\def\genspazio #1#2#3#4#5{#1^{#2}(#5,#4;#3)}
\def\spazio #1#2#3{\genspazio {#1}{#2}{#3}T0}

\def\L {\spazio L}
\def\H {\spazio H}
\def\W {\spazio W}

\def\C #1#2{C^{#1}([0,T];#2)}

\def\Vp{V^*}


\def\Lx #1{L^{#1}(\Omega)}
\def\Hx #1{H^{#1}(\Omega)}

\def\Cx #1{C^{#1}(\overline\Omega)}

\def\Ldue{\Lx 2}
\def\Linfty{\Lx\infty}

\def\Huno{\Hx 1}
\def\Hdue{\Hx 2}


\def\LQ #1{L^{#1}(Q)}


\let\theta\vartheta
\let\eps\varepsilon

\let\TeXchi\chi                         
\newbox\chibox
\setbox0 \hbox{\mathsurround0pt $\TeXchi$}
\setbox\chibox \hbox{\raise\dp0 \box 0 }
\def\chi{\copy\chibox}


\def\muz{\mu_0}
\def\rhoz{\rho_0}

\def\rhomin{\rho_*}
\def\rhomax{\rho^*}

\def\normaV #1{\norma{#1}_V}
\def\normaH #1{\norma{#1}_H}

\def\kmin{\kappa_*}

\def\coeff{1+2g(\rho)}

\def\coefft{1+2g(\rho(t))}

\def\mus{\mu_\sigma}
\def\rhos{\rho_\sigma}
\def\xis{\xi_\sigma}

\def\an{a_n}
\def\bn{b_n}
\def\zn{z_n}




\def\Krejci{Krej\v c\'\i{}}

\Begin{document}


\title{\bf Regularity of the solution 
to a nonstandard system of phase field equations%
\footnote{{\bf Acknowledgments.}\quad\rm 
PC and GG gratefully acknowledge some financial support from the MIUR-PRIN Grant 2010A2TFX2 ``Calculus of Variations'' and the GNAMPA (Gruppo Nazionale per l'Analisi Matematica, la Probabilit\`a e le loro Applicazioni) of INdAM (Istituto Nazionale di 
Alta Matematica).}}
\author{}
\date{}
\maketitle
\begin{center}
\vskip-2cm
{\pier
{\large\bf Pierluigi Colli$^{(1)}$}\\
{\normalsize e-mail: {\tt pierluigi.colli@unipv.it}}\\[.25cm]
{\large\bf Gianni Gilardi$^{(1)}$}\\
{\normalsize e-mail: {\tt gianni.gilardi@unipv.it}}\\[.25cm]
{\large\bf J\"urgen Sprekels$^{(2)}$}\\
{\normalsize e-mail: {\tt sprekels@wias-berlin.de}}\\[.45cm]
$^{(1)}$
{\small Dipartimento di Matematica ``F. Casorati'', Universit\`a di Pavia}\\
{\small via Ferrata 1, 27100 Pavia, Italy}\\[.2cm]
$^{(2)}$
{\small Weierstra\ss-Institut f\"ur Angewandte Analysis und Stochastik}\\
{\small Mohrenstra\ss e\ 39, 10117 Berlin, Germany}\\[.8cm]
}
\end{center}


\Begin{abstract}
A nonstandard system of differential equations
describing two-species phase segregation is considered.
{\juerg This} system naturally arises
in the asymptotic analysis recently done by
Colli, Gilardi, \Krejci, and Sprekels as the diffusion coefficient 
in the equation governing the evolution of the order parameter tends to zero.
In particular, a well-posedness result is proved for the limit system.
This paper deals with the above limit problem 
in a less general but still very significant framework
and provides a very simple proof of further regularity for the solution.
As a byproduct, a~simple uniqueness proof is given as well.

\medskip\noindent
{\bf Key words:}
nonstandard phase field system, nonlinear differential equations, 
{\pier initial and boundary value problem, regularity of solutions.}\\[0.2cm]
{\bf AMS (MOS) Subject Classification:} 35K61, 35D10, 35A02.
\End{abstract}


\salta

\pagestyle{myheadings}
\newcommand\testopari{\sc Colli \ --- \ Gilardi \ --- \ Sprekels}
\newcommand\testodispari{\sc Regularity for a nonstandard phase field {\juerg system}}
\markboth{\testodispari}{\testopari}

\finqui


\section{Introduction}
\label{Intro}
\setcounter{equation}{0}
In this paper, we consider the {\juerg system} 
\Bsist
  && \bigl( 1 + 2g(\rho) \bigr) \, \dt\mu
  + \mu \, g'(\rho) \, \dt\rho
  - \Delta\mu = 0
  \label{Iprima}
  \\
  && \dt\rho + f'(\rho) = \mu \, g'(\rho)
  \label{Iseconda}
  \\
  && \dn\mu|_\Gamma = 0
  \label{Ibc}
  \\
  && \mu(0) = \muz
  \aand
  \rho(0) = \rhoz 
  \label{Icauchy}
\Esist
\Accorpa\Ipbl Iprima Icauchy
in the unknown fields $\mu$ and~$\rho$,
where {\juerg the} equations \accorpa{Iprima}{Iseconda} are meant to hold in a 
bounded domain $\Omega\subset\erre^3$ with a smooth boundary~$\Gamma$ 
and in some time interval~$(0,T)$,
and {\juerg where} $\dn$ in \eqref{Ibc} denotes the {\juerg outward} normal derivative. 
In the recent {\pier papers~\cite{CGKS1,CGKS2}, the well-posedness of 
the problem \eqref{Iprima}--\eqref{Icauchy} was investigated, and in particular the existence of the solution} {\juerg was} proved
by considering the {\juerg system} of partial differential equations
obtained by replacing the ordinary differential equation~\eqref{Iseconda}
by the partial differential equation
\Beq
  \dt\rho - \sigma \Delta\rho + f'(\rho) = \mu \, g'(\rho)
  \quad \hbox{with the boundary condition} \quad
  \dn\rho|_\Gamma = 0
  \label{Isecondas}
\Eeq
and performing the asymptotic analysis as $\sigma$ tends to zero.
{\juerg This} modified system originates from the model introduced in~\cite{Podio},
{\juerg which} describes the phase segregation of two species 
(atoms and vacancies, say) 
on~a lattice in presence of diffusion
and looks like a modification of the well-known Cahn-Hilliard equations
(see, e.g.,~{\pier \cite{FG, Gurtin}}).
The~state variables are the  {\sl order parameter\/}~$\rho$
(volume density of one of the two species),
which {\juerg of course must attain} values in the domain of~$f'$,
and the {\sl chemical potential\/}~$\mu$,
which is required to be nonnegative for physical reasons. 
{\juerg This} system has been studied in a series of papers 
and a number of results has been obtained
in several directions {\pier \cite{CGKPS, CGPS3, CGPS4, CGPS5, CGPS6, CGS1}.}
Moreover, some of {\juerg these} results hold for a more general system 
involving a nonlinear elliptic operator in divergence form
in equation~\eqref{Iprima}, {\juerg in place} of the Laplacian {\pier (see \cite{CGKS1, CGKS2, CGPS7, CGPS8})}.
In all {\juerg of} the quoted papers, the function $f$ represents
a double-well potential{\juerg . A} thermodynamically relevant example is
the so-called {\sl logarithmic potential\/} 
defined (up~to an additive constant) by~the formula
\Beq
  f(\rho) = c_1 \, \bigl( \rho\,\log\rho + (1-\rho) \, \log(1-\rho) \bigl)
  {}+{} c_2 \, \rho (1-\rho)
  \quad \hbox{for $\rho\in(0,1)$},
  \label{logpot} 
\Eeq
where $c_1$ and $c_2$ are positive constant
with $c_2>2c_1$
in order that $f$ actually presents a double well.
However, the class of the admissible potentials 
could be much wider and include{\juerg s} both the standard double-well potential defined~by
\Beq
  f(\rho) = \frac 14 \, (\rho^2-1)^2
  \quad \hbox{for $\rho\in\erre$}
  \label{standardpot}
\Eeq
and potentials whose convex part is just proper and lower semicontinuous, 
{\juerg and} thus possibly non-differentiable, in its effective domain.
In the latter case,
the~monotone part of $f'$ is replaced by a multivalued subdifferential 
and \eqref{Isecondas} has to be read as a differential inclusion. 
In~\cite{CGKS1}, such a wide class of potentials is considered,
so that the existence result for system \Ipbl\ obtained there is very general.
However, the solution constructed in this way may be irregular, in principle.
Nevertheless, it is unique and a little more regular than expected,
{\juerg at} the price that the corresponding proofs are rather complicated.
 
The present paper deals just with potentials
that see example \eqref{logpot} as a prototype,
but it provides simple proofs of further regularity.
As an application, we give an easy uniqueness proof.

Our paper is \organiz ed as follows.
In the next section, we list our assumptions 
and state problem \Ipbl\ in a precise form.
In the last section, we present and prove our results.


\section{Assumptions and notations}
\label{MainResults}
\setcounter{equation}{0}

We first introduce precise assumptions on the data 
for the mathematical problem under investigation. 
We assume $\Omega$ to be a bounded connected 
open set in $\erre^3$ with smooth boundary~$\Gamma$
({treating} lower-dimensional cases would require {only }minor changes) 
and let  $T\in(0,+\infty)$ stand for a final time. 
We~set 
\Beq
  V := \Huno, \quad
  H := \Ldue , 
  \aand
  W := \graffe{v\in\Hdue:\ \dn v|_\Gamma = 0},
  \label{defspazi}
\Eeq
and endow the spaces \eqref{defspazi} with their standard norms,
for which we use a self-explanato\-ry notation like $\normaV\cpto$.
For simplicity, we use the same {\juerg notation also for powers of these spaces}.
The symbol $\<\cpto,\cpto>$ denotes the duality product 
between~$\Vp$, the dual space of~$V$, and~$V$ itself.
Moreover, for $p\in[1,+\infty]$, we write $\norma\cpto_p$ 
for the usual norm both in $\Lx p$ and in~$\LQ p$,
where $Q:=\Omega\times(0,T)$.
For the nonlinearities and the initial data we assume that
there exist
\Beq
  \rhomin , \rhomax \in \erre \quad \hbox{with} \quad
  \rhomin < \rhomax
  \label{hprhomm}
\Eeq
in order that the combined conditions listed below hold.
\Bsist
  && f ,\, g : [\rhomin,\rhomax] \to \erre 
  \quad \hbox{are $C^2$ functions}
  \label{hpfg}
  \\
  && g(r) \geq 0 
  \aand
  g''(r) \leq 0
  \quad \hbox{for every $r\in[\rhomin,\rhomax]$} 
  \qquad
  \label{hpg}
  \\
  && f'(\rhomin) \leq 0 \leq g'(\rhomin)
  \aand
  g'(\rhomax) \leq 0 \leq f'(\rhomax) .
  \label{hpsign}
  \\
  && \muz \in V \cap \Linfty
  \aand
  \muz \geq 0
  \quad \aeO
  \label{hpmuz}
  \\
  && \rhoz \in V 
  \aand
  \rhomin \leq \rhoz \leq \rhomax
  \quad \aeO \quad
  \label{hprhoz}
\Esist
\Accorpa\Hptutto hprhomm hprhoz
Notice that the functions $f$, $g$, {\juerg together with their first derivatives,}
are bounded and \Lip\ continuous.
Moreover, we remark that the different assumptions of \cite{CGKS1} 
can be fulfilled by splitting $f$ as $f_1+f_2$
with $f_1$ nonnegative and convex, and suitably extending $f_1$, $f_2$, and $g$ 
to an open interval including~$[\rhomin,\rhomax]$.
In particular, the logarithmic potential \eqref{logpot} 
{\juerg fits} the above requirements with $\rhomin,\rhomax\in(0,1)$ and reasonable choices of~$g$.

Now, we recall the part that follows from the asymptotic analysis performed in~\cite{CGKS1}
and is of interest for the present paper.

\Bthm
\label{Esistenza}
Assume that {\juerg the assumptions}~\Hptutto\ hold. 
Then there exists at least {\juerg one} pair $(\mu,\rho)$ satisfying
\Bsist 
  && \mu \in \L\infty H \cap \L2V \cap \LQ\infty
  \aand
  \mu \geq 0 \quad \aeQ
  \label{regmu}
  \\
  && \rho \in \H1H \cap \L\infty V
  \aand
  \rhomin \leq \rho \leq \rhomax \quad \aeQ
  \qquad
  \label{regrho}
  \\
  && u := \bigl( 1+2g(\rho) \bigr) \mu \in \W{1,1}\Vp
  \label{regu}
\Esist
\Accorpa\Regsoluz {regmu} {regu}
and solving the problem 
\Bsist
  && \< \dt u(t), v >
  + \int_\Omega  \nabla\mu(t) \cdot \nabla v
  = \iO \mu(t) \, g'(\rho(t)) \, \dt\rho(t) \, v 
  \non
  \\
  && \qquad
  \hbox{for all $v\in V$ and a.a.\ $t\in(0,T)$}
  \label{prima}
  \\
  && \dt\rho + f'(\rho)
   = \mu \, g'(\rho)
  \quad \aeQ
  \label{seconda}
  \\
  && u(0) = \bigl( 1+2g(\rhoz) \bigr) \muz
  \aand
  \rho(0) = \rhoz
  \quad \aeO .
  \label{cauchy}
\Esist
\Accorpa\Pbl prima cauchy
\Ethm

We observe that the first regularity level for the time derivative of $u$
obtained in \cite{CGKS1} is $\dt u\in\L{4/3}\Vp$, that is,
a~little better than~\eqref{regu}.
However, one easily sees that $\dt u\in\L2\Vp$
by comparison in~\eqref{prima},
on account of~\accorpa{regmu}{regrho} and~\eqref{hpfg}.

\Brem
\label{Primastrong}
We can also consider the stronger form of~\eqref{prima},
\Beq
  \iO \bigl( \coefft \bigr) \, \dt\mu(t) \, v
  + \iO \mu (t)\, g'(\rho(t)) \, \dt\rho(t) \, v   
  + \int_\Omega  \nabla\mu(t) \cdot \nabla v
  = 0
  \label{primastrong}
\Eeq
for all $v\in V$ and \aat,
and observe that it is equivalent to~\eqref{prima}
provided that one can apply the Leibniz rule to the time derivative~$\dt u$.
This is the case if $\dt\mu$ exists and belongs to~$\LQ2$.
However, Theorem~\ref{Esistenza} does not ensure such a regularity.
Moreover, \eqref{prima}~also includes the homogenous Neumann boundary condition~\eqref{Ibc}
in a \generaliz ed sense.
\Erem

The aim of this paper is to prove that any solution to problem~\Pbl\ 
satisfying the very mild regularity~\Regsoluz\
is in fact much smoother and, in particular, unique.

Now, we list a number of tools and notations {\juerg used} throughout the paper.
First of all, we often use the elementary Young inequality,
\Beq
  ab\leq \eps a^2 + \frac 1 {4\eps} \, b^2
  \quad \hbox{for every $a,b\geq 0$ and $\eps>0$},
  \label{young}
\Eeq
and repeatedly account for the H\"older and Sobolev inequalities.
The precise form of the latter we use is the following:
\Beq
  V \subset \Lx q
  \aand
  \norma v_{q} \leq C \normaV v
  \quad \hbox{for every $v\in V$ and $q\in[1,6]$},
  \label{sobolev}
\Eeq
where $C$~depends only on~$\Omega$.
Moreover, the above embedding is compact if $q<6$,
and the compactness inequality
\Beq
  \norma v_q \leq \eps \norma{\nabla v}_2 + C_{q,\eps} \norma v_2
  \quad \hbox{for every $v\in V$, $q\in[1,6)$, and $\eps>0$},
  \label{compact}
\Eeq
holds with a constant $C_{q,\eps}$ depending on $\Omega$, $q$, and~$\eps$, only.
Furthermore, we exploit the embeddings
\Beq
  \L\infty H \cap \L2V
  \subset \L\infty H \cap \L2{\Lx6}
  \subset \LQ{10/3},
  \label{incldieciterzi}
\Eeq 
{\juerg as well as} the corresponding inequality
\Beq
  \norma v_{\LQ{10/3}}
  \leq C \bigl(
    \norma v_{\L\infty H} + \norma v_{\L2V}
  \bigr),
  \label{disdieciterzi}  
\Eeq
{\juerg which follow from} combining the Sobolev embedding $V\subset\Lx6$
and the well-known interpolation inequality for $L^p$ spaces.
Again, $C$~depends only on~$\Omega$.
Finally, we set 
\Beq
  Q_t := \Omega \times (0,t)
  \quad \hbox{for $t\in(0,T]$}
  \label{defQt}
\Eeq
and use the same symbol small-case $c$ for different constants, 
that may only depend on~$\Omega$, the final time~$T$, the nonlinearities $f$ and~$g$, 
and the properties of the data involved in the statements at hand.
A~notation like~$c_\eps$ signals a constant that {\juerg also depends} on the parameter~$\eps$. 
The reader should keep in mind that the meaning of $c$ and $c_\eps$ {\juerg may}
change from line to line and even {\juerg within} the same chain of inequalities, 
whereas those constants we need to refer to are always denoted by 
capital letters, just like $C$ in~\eqref{sobolev} and in~\eqref{disdieciterzi}.


\section{Regularity}
\label{Results}
\setcounter{equation}{0}

In this section, we prove regularity results for the solution to problem \Pbl\
{\juerg under the assumption that the} conditions \Hptutto\ hold (we~often avoid writing this).
In~order to help the reader, we sketch our strategy.
We fix any solution $(\mu,\rho)$ to problem \Pbl\ satisfying the regularity requirements \Regsoluz\
and recall that all {\juerg of} the nonlinear terms involving $\rho$ are bounded.
Moreover, $\mu$~is bounded, too (cf.~\eqref{regmu}).
Thus, \eqref{seconda}~implies that even $\dt\rho$ is bounded.
Now, we~set
\Beq
  a := \coeff
  \aand
  b := \mu \, g'(\rho) \, \dt\rho 
  \label{defab}
\Eeq
and notice that $\dt a=2g'(\rho)\dt\rho$.
Hence, we have
\Beq
  a \in \L\infty V \cap \LQ\infty, \quad
  \dt a \in \LQ\infty , \quad
  b \in \LQ\infty , \quad 
  a \geq 1 \quad \aeQ .
  \label{regab}
\Eeq
Next, we introduce the associated linear problem
\Bsist
  && \< \dt (az)(t) , v >
  + \iO  \nabla z(t) \cdot \nabla v
  = \iO b(t) \, v 
  \quad \hbox{for all $v\in V$ and a.a.\ $t\in(0,T)$}
  \qquad
  \label{primaz}
  \\
  && (az)(0) = \bigl( 1+2g(\rhoz) \bigr) \muz \,, 
  \label{cauchyz}
\Esist
\Accorpa\Pblz {primaz} {cauchyz}
whose unknown $z$ is required to satisfy
\Beq
  z \in \L\infty H \cap \L2V 
  \aand
  az \in \W{1,1}\Vp,
  \label{regz}
\Eeq
and observe that $z=\mu$ is a solution.
Then, we prove that \Pblz\ has a unique solution $z$
satisfying~\eqref{regz}.
This implies the following.
If we \regulariz e \Pblz\ and perform some 
a~priori estimates on the solution to the \regulariz ed problem {\juerg then
these} estimates still hold for any weak limit.
On the other hand, such a limit must be~$\mu$ due to uniqueness.
This {\juerg entails further} regularity for~$\mu$.
Once the regularity obtained for $\mu$ is sufficiently high,
we can even prove uniqueness in a simple way.
We observe that uniqueness for \Pblz\ is not \sfw.
Indeed, \eqref{primaz}~is a very weak form 
(due~to the very low regularity~\eqref{regz})
of~the homogeneous Neumann boundary value problem 
for the equation
\Beq
  \dt(az) - \Delta z = b\,,
  \non
\Eeq
which is formally uniformly parabolic.
However, the equation is not presented in divergence form,
and $a$ might be discontinuous
since no continuity for $\rho$ is known.
At this point, we can start with our program.

\Blem
\label{Uniclineare}
Let $(\mu,\rho)$ be a solution to \Pbl\ satisfying \Regsoluz,
and let $a$ and $b$ be defined by~\eqref{defab}.
Then problem \Pblz\ has a unique solution $z$ satisfying~\eqref{regz},
and {\juerg this solution coincides with}~$\mu$.
\Elem

\Bdim
Clearly, $\mu$~satisfies \eqref{regz} and solves~\Pblz.
As far as uniqueness is concerned, 
we can deal with the correponding homogeneous problem, by linearity.
Thus, we fix any $z$ satisfying \eqref{regz} that solves~\Pblz, 
where $b$ is replaced by {\juerg zero} on the \rhs\ of \eqref{primaz}
and the initial value in \eqref{cauchyz} is~{\juerg zero}.
Then, we introduce the adjoint problem of finding $v$ such~that
\Bsist
  && v \in \H1H \cap \L2W \subset \C0V
  \label{regv}
  \\
  && - a \, \dt v - \Delta v = z \quad \aeQ
  \aand
  v(T) = 0 .
  \label{aggiunto}
\Esist
It is easily seen that \accorpa{regv}{aggiunto} has a (unique) solution
since $z\in\L2H$.
We use such a solution $v$ as a test function for~$z${\juerg , 
observing} that the integration by parts {\juerg formula}
\Beq
  \ioT \< \dt w(t) , \phi(t) > \, dt
  = \< w(T) , \phi(T) > 
  - \< w(0) , \phi(0) > 
  - \ioT \< \dt\phi(t) , w(t) > \, dt
  \non
\Eeq
is actually {\juerg valid} if $w\in\L2V\cap\W{1,1}\Vp$ and $\phi\in\C0V\cap\H1\Vp$.
{\juerg Thus}, the choice $\phi=v$ is allowed by~\eqref{regv},
Hence, we~have
\Bsist
  && 0 = \ioT \< \dt(az)(t) , v(t) > \, dt
  + \intQ \nabla z \cdot \nabla v
  \non
  \\
  && = \< (az)(T) , v(T) > 
  - \< (az)(0) , v(0) > 
  - \ioT \< \dt v(t) , (az)(t) > \, dt
  - \intQ z \, \Delta v
  \non
  \\
  && = \intQ \Bigl( -a \, \dt v - \Delta v \Bigr) z
  = \intQ z^2\,,
  \non
\Esist
{\juerg which implies that} $z=0$.
\Edim

\Bcor
\label{Onlymu}
If $z\in\H1H\cap\L\infty V\cap\L2W$ solves
\Beq
  a \, \dt z + \dt a \, z - \Delta z = b
  \quad \aeQ
  \aand z(0) = \muz,
  \label{pblzforte}
\Eeq
then {\juerg it holds} $z=\mu$.
\Ecor

\Bdim
Indeed, our assumptions on $z$ and \eqref{pblzforte} 
imply both \eqref{regz} and \Pblz. 
Thus, the previous lemma gives the result.
\Edim

\Bthm
\label{Regolarita}
Assume that {\juerg the assumptions}~\Hptutto\ {\juerg are fulfilled,}
and let $(\mu,\rho)$ be a solution to \Pbl\ satisfying \Regsoluz.
Then $\mu$ also satisfies
\Beq
  \mu \in \H1H \cap \L2W .
  \label{regolarita}
\Eeq
In particular, $\mu\in\C0V$.
\Ethm

\Bdim
Thanks to \eqref{regab},
there exist two sequences $\graffe{a_n}$ and $\graffe{b_n}$ of $C^1$ functions
such~that
\Bsist
  && \an \to a , \quad \dt\an \to \dt a , \quad \bn \to b
  \quad \hbox{strongly in $\LQ2$}
  \label{convn} 
  \\
  && |\an| + |\dt\an| + |\bn| \leq c
  \aand
  \an \geq 1 \quad \aeQ
  \label{limitn}
\Esist
Then, we consider {\juerg for any $n\in\enne$} the following \regulariz ation of problem~\eqref{pblzforte}:
\Beq
  \an \dt\zn + \dt\an \, \zn - \Delta\zn = \bn 
  \quad \aeQ
  \aand
  \zn(0) = \muz \,,
  \label{pblzn}
\Eeq
complemented with the homogeneous Neumann boundary condition.
Due to the regularity of the coefficients and assumption \eqref{hpmuz} on~$\muz$,
it is easy to see that {\juerg this} problem has~a (unique) solution~$\zn$
satisfying 
\Beq
  \zn \in \H1H \cap \L\infty V \cap \L2W.
  \label{regzn}
\Eeq
We add $\zn$ to both sides of the equation for convenience.
Then, we multiply the resulting equality by~$\dt\zn$ 
and integrate over~$Q_t$ {\juerg for} any $t\in(0,T)$.
Owing to~\eqref{limitn}, we easily obtain {\juerg that}
\Bsist
  && \intQt |\dt\zn|^2
  + \frac 12 \normaV{\zn(t)}^2
  \leq \frac 12 \iO |\nabla\muz|^2
  + \intQt \bigl( \bn + \zn - \dt\an \, \zn) \dt\zn
  \non
  \\
  && \leq \frac 12 \iO |\nabla\muz|^2
  + \frac 12 \intQt |\dt\zn|^2
  + c \intQt (1+|\zn|^2)
  \non
  \\
  && \leq c
  + \frac 12 \intQt |\dt\zn|^2
  + c \iot \normaV{\zn(s)}^2 \, ds .
  \non
\Esist
By rearranging and applying the Gronwall lemma,
we immediately conclude that
\Beq
  \norma{\dt\zn}_{\L2H}
  + \norma\zn_{\L\infty V}
  \leq c .
  \non
\Eeq
Moreover, by comparison in~\eqref{pblzn},
we also find an estimate for~$\normaH{\Delta\zn}$.
So, the above bounds and standard regularity results for elliptic equations yield
{\juerg the estimate}
\Beq
  \norma\zn_{\H1H}
  + \norma\zn_{\L\infty V}
  + \norma\zn_{\L2W}
  \leq c .
  \label{primastima}
\Eeq
Therefore, by weak compactness, there exists {\juerg some} $z$ such that
\Beq
  \zn \to z
  \quad \hbox{weakly star in $\H1H\cap\L\infty V\cap\L2W$},
  \non
\Eeq
at least for a subsequence.
This also implies $\zn\to z$ weakly in~$\C0H$
and, recalling~\eqref{convn}, we infer that
\Beq
  a \, \dt z + \dt a \, z - \Delta z = b
  \aand
  z(0) = \muz \,.
  \non
\Eeq
Now, we apply Corollary~\ref{Onlymu} and conclude that $z=\mu$, whence~\eqref{regolarita}
{\juerg follows}.
The last {\juerg assertion {\pier is a consequence of}} the well-known embedding 
$\H1H\cap\L2W\subset\C0V.$
\Edim

The regularity just established can be improved
provided that a stronger assumption on the initial datum~$\muz$ is satisfied, namely
\Beq
  \muz \in W .
  \label{hpW}
\Eeq
We observe that the regularity given in \eqref{hpW} implies \eqref{hpmuz} 
because of the continuous embedding $W\subset\Cx0$.
The new regularity result is stated in the following theorem,
whose proof is performed with the same technique as before.

\Bthm
\label{Piuregolarita}
In addition to~\Hptutto, assume that \eqref{hpW} holds,
and let $(\mu,\rho)$ be a solution to \Pbl\ satisfying \Regsoluz.
Then $\mu$ also satisfies
\Beq
  \dt\mu \in \L\infty H \cap \L2V 
  \aand
  \mu \in \L\infty W .
  \label{piuregolarita}
\Eeq
\Ethm

\Bdim
Theorem~\ref{Regolarita} ensures that \eqref{regolarita}~holds for~$\mu$,
so that the regularity of $\dt\rho$, $a$, and~$b$ can be updated.
Indeed, \eqref{regolarita}~implies that \eqref{seconda} can be differentiated with respect to time
and that $\dt^2\rho\in\LQ2$.
Hence, we also~have
\Beq
  \dt^2 a \in \LQ2
  \aand 
  \dt b \in \LQ2, 
  \label{piuregab}
\Eeq
{\juerg which} allows us to construct a new approximation of problem~\eqref{pblzforte}.
Namely, we can choose sequences $\graffe{a_n}$ and $\graffe{b_n}$ of $C^2$ functions
satisfying \accorpa{convn}{limitn}~and
\Beq
  \norma{\dt^2\an}_{\LQ2} + \norma{\dt\bn}_{\LQ2}
  \leq c
  \quad \hbox{for every $n$} .
  \label{piulimitn}
\Eeq
This leads to a sequence of solutions $\zn$ to the corresponding problems~\eqref{pblzn},
which, however, still keeps all {\juerg of} the properties of the approximation 
we had established in the proof of the previous theorem.
But the regularity of the coefficient and {\juerg the} assumption \eqref{hpW} on $\muz$
also {\juerg ensure} us that $\zn$ is smoother 
and that equation~\eqref{pblzn} can be differentiated with respect to time.
By doing this, we~get
\Beq
  \an \, \dt^2\zn + 2 \dt\an \, \dt\zn + \dt^2\an \, \zn - \Delta\dt\zn = \dt\bn 
  \non\,,
\Eeq
and we can multiply {\juerg this} equality by~$\dt\zn$ and integrate over~$Q_t$,
where $t\in(0,T)$ is arbitrary.
By integrating by parts with respect to time 
the term involving the second derivative~$\dt^2\zn$,
and owing to the inequality $\an\geq1$,
we obtain
\Bsist
  && \frac 12 \iO |\dt\zn(t)|^2 
  + \intQt |\nabla\dt\zn|^2
  \non
  \\
  && \leq \frac 12 \iO \an(0) |\dt\zn(0)|^2 
  - \frac 52 \intQt \dt\an \, |\dt\zn|^2
  - \intQt \dt^2\an \, \zn \, \dt\zn
  + \intQt \dt\bn \, \dt\zn \,.
  \qquad\quad
  \label{persecondastima}
\Esist
Two of the terms on the \rhs\ can be dealt with{\juerg , 
by accounting for \eqref{limitn}, \eqref{primastima} and \eqref{piulimitn}, in the following way:}
\Beq
  - \frac 52 \intQt \dt\an \, |\dt\zn|^2
  + \intQt \dt\bn \, \dt\zn 
  \leq  c + c \intQt |\dt\zn|^2
  \leq c .
  \non
\Eeq
For the first one, we observe that \eqref{pblzn}, \eqref{hpW}, and \eqref{limitn} imply
{\juerg that} $\an(0)\leq c$ and
\Beq
  |\dt\zn(0)|
  \leq |\bn(0) + \Delta\zn(0) - \dt\an(0) \, \zn(0)|
  \leq c + |\Delta\muz| + c \, \muz
  \leq c + |\Delta\muz| .
  \non
\Eeq
As $\muz\in W$, the integral under investigation is bounded.
It remains to handle the third term.
By using first the \holder\ and Young inequalities,
and then the Sobolev and compactness inequalities~\eqref{sobolev} and \eqref{compact} with $q=4$, 
we~have 
\Bsist
  && - \intQt \dt^2\an \, \zn \, \dt\zn
  \leq \iot \norma{\dt^2\an(s)}_2 \, \norma{\zn(s)}_4 \, \norma{\dt\zn(s)}_4 \, ds
  \non
  \\
  && \leq \iot \norma{\dt\zn(s)}_4^2
  + \iot \norma{\dt^2\an(s)}_2^2 \, \norma{\zn(s)}_4^2 \, ds
  \non
  \\
  && \leq \frac 12 \iot \normaH{\nabla\dt\zn(s)}^2 \, ds
  + c \iot \normaH{\dt\zn(s)}^2 \, ds
  + c \iot \norma{\dt^2\an(s)}_2^2 \, \normaV{\zn(s)}^2 \, ds 
  \non
  \\
  && \leq \frac 12 \iot \normaH{\nabla\dt\zn(s)}^2 \, ds
  + c\,,
  \non
\Esist
{\juerg where the last inequality follows from}~\eqref{primastima} and~\eqref{piulimitn}.
By collecting all this and \eqref{persecondastima} and rearranging,
we~conclude~that
\Beq
  \norma{\dt\zn}_{\L\infty H} + \norma{\dt\zn}_{\L2V} \leq c .
  \label{secondastima}
\Eeq
Moreover, by comparison in~\eqref{pblzn},
we also infer that {\juerg $\{\Delta\zn\}$} is bounded in~$\L\infty H$.
Elliptic regularity and \eqref{primastima} then yield {\juerg the boundedness of 
$\{\zn\}$} in~$\L\infty W$.
At this point, we use weak compactness once more.
We obtain{\juerg , on a subsequence, that} 
\Beq
  \zn \to z
  \quad \hbox{weakly star in $\W{1,\infty}H\cap\H1V\cap\L\infty W$}
  \non
\Eeq 
and conclude that $z=\mu$ as before.
Hence, \eqref{regolarita} holds, and the proof is complete.
\Edim

\Bcor
\label{Dieciterzi}
Under the assumptions of Theorem~\ref{Piuregolarita},
every solution $(\mu,\rho)$ to problem \Pbl\ also satisfies
\Beq
  \dt\mu \in \LQ{10/3} .
  \label{dieciterzi}
\Eeq
\Ecor

\Bdim
It suffices to combine the first of \eqref{piuregolarita} and \eqref{incldieciterzi}.
\Edim

The regularity just achieved allows us to give a rather simple uniqueness proof.

\Bthm
\label{Uniqueness}
Assume that~\Hptutto\ and \eqref{hpW} hold. 
Then the solution $(\mu,\rho)$ given by Theorem~\ref{Esistenza} is unique.
\Ethm

\Bdim
We adapt the argument used in~\cite{CGPS3} to the present situation.
We pick two solutions $(\mu_i,\rho_i)$, $i=1,2$,
and set for convenience
\Beq
  \rho := \rho_1-\rho_2 \,, \enskip
  \mu := \mu_1-\mu_2 \,, \enskip
  \gamma_i := g(\rho_i) \,, \enskip
  \gamma := \gamma_1-\gamma_2 \,, \enskip
  \eta_i := g'(\rho_i) \,, \enskip
  \eta := \eta_1-\eta_2 \,.
  \non
\Eeq
By accounting for Remark~\ref{Primastrong},
we write \eqref{primastrong} for both solutions, take the difference, and multiply by~$\mu$.
At the same time, we write \eqref{seconda} for both solutions,
add $\rho_i$ to both sides for convenience,
move all the nonlinear terms to the \rhs,
and multiply the difference by~$\dt\rho$.
Then, we integrate and sum up.
By owing to the identity
\Bsist
  && \{
    (1+2\gamma_1) \dt\mu_1 + \mu_1 \eta_1 \dt\rho_1
    - (1+2\gamma_2) \dt\mu_2 + \mu_2 \eta_2 \dt\rho_2
  \} \, \mu
  \non
  \\
  && = \frac 12 \, \dt \{ (1+2\gamma_1) |\mu|^2 \}
  + 2\gamma \dt\mu_2 \, \mu
  + \mu_2 \eta \dt\rho_1 \, \mu
  + \mu_2 \eta_2 \dt\rho \, \mu\,,
  \label{checkplease}
\Esist
and using the boundedness and the \Lip\ continuity of the nonlinearities,
we~obtain
\Bsist
  && \hskip -4pt 
  \frac 12 \iO |\mu(t)|^2
  + \intQt |\nabla\mu|^2 
  + \intQt |\dt\rho|^2
  + \frac 12 \iO |\rho(t)|^2
  \non
  \\
  && \leq c \intQt |\dt\mu_2| \, |\rho| \, |\mu|
   + c \intQt \bigl(
   |\rho| \, |\mu|
   + |\dt\rho| \, |\mu|
   + |\dt\rho| \, |\rho|
  \bigr) .
  \label{perunicita}
\Esist
Just the first integral on the \rhs\ needs some treatement.
By the \holder\ inequality, the compactness inequality \eqref{compact} with $q=5$,
and the elementary Young inequality, 
we~have for any~$\eps>0$
\Bsist
  && \intQt |\dt\mu_2| \, |\rho| \, |\mu|
  \leq \iot \norma{\dt\mu_2(s)}_{10/3} \, \norma{\rho(s)}_2 \, \norma{\mu(s)}_5 \, ds
  \non
  \\
  && \leq \iot \norma{\mu(s)}_5^2 \, ds
  + \iot \norma{\dt\mu_2(s)}_{10/3}^2 \, \norma{\rho(s)}_2^2 \, ds 
  \non
  \\
  && \leq \eps \iot \norma{\nabla\mu(s)}_2^2 \, ds
  + c_\eps \iot \norma{\mu(s)}_2^2 \, ds
  + \iot \norma{\dt\mu_2(s)}_{10/3}^2 \, \norma{\rho(s)}_2^2 \, ds .
  \non
\Esist
As the function $s\mapsto\norma{\dt\mu_2(s)}_{10/3}^2$ belongs 
to $L^{5/3}(0,T)$ by~Corollary~\ref{Dieciterzi}, {\juerg and} thus to~$L^1(0,T)$,
the last integral can be controlled by the \lhs\ of \eqref{perunicita}
via Gronwall-Bellman's lemma (see, e.g., \cite[Lemma~A.4, p.~156]{Brezis}).
Hence, it is sufficient to choose $\eps$ small enough and apply {\juerg this} lemma
in order to conclude that $\mu=0$ and $\rho=0$.
\Edim

\Brem
It is clear that the bootstrap procedure used in the above proofs
can {\juerg be continued to provide} even more regularity for the solution $(\mu,\rho)$ to problem \Pbl\
under suitable assumptions on the initial data.
\Erem



\vspace{3truemm}

\Begin{thebibliography}{10}

\bibitem{Brezis}
H. Brezis,
``Op\'erateurs maximaux monotones et semi-groupes de contractions
dans les espaces de Hilbert'',
North-Holland Math. Stud.
{\bf 5},
North-Holland,
Amsterdam,
1973.

{\pier
\bibitem{CGKPS} 
P. Colli, G. Gilardi, P. \Krejci, P. Podio-Guidugli, J. Sprekels,
Analysis of a time discretization scheme for
a nonstandard viscous Cahn-Hilliard system, 
{\it ESAIM Math. Model. Numer. Anal.} {\bf 48} (2014) 1061-1087.}

\bibitem{CGKS1} 
P. Colli, G. Gilardi, P. \Krejci, J. Sprekels,
A vanishing diffusion limit 
in a nonstandard system of phase field equations,
{\it Evol. Equ. Control Theory\/} {\bf 3} (2014) 257-275.

{\pier
\bibitem{CGKS2} 
P. Colli, G. Gilardi, P. \Krejci, J. Sprekels,
A continuous dependence result
for a nonstandard system of phase field equations,
{\it Math. Methods Appl. Sci.} {\bf 37} (2014) 1318-1324.}

\bibitem{CGPS3} 
P. Colli, G. Gilardi, P. Podio-Guidugli, J. Sprekels,
Well-posedness and long-time behaviour for a 
nonstandard viscous Cahn-Hilliard system, 
{\it SIAM J. Appl. Math.} {\bf 71} (2011) 1849-1870.

{\pier
\bibitem{CGPS4} 
{P. Colli, G. Gilardi, P. Podio-Guidugli, J. Sprekels},
An asymptotic analysis for a nonstandard Cahn-Hilliard system with viscosity, 
{\it  Discrete Contin. Dyn. Syst. Ser. S}  {\bf 6} (2013) 353-368. 

\bibitem{CGPS5} 
{P. Colli, G. Gilardi, P. Podio-Guidugli, J. Sprekels},
Distributed optimal control of a nonstandard 
system of phase field equations, 
{\it Contin. Mech. Thermodyn.} {\bf 24} (2012) 437-459.
}

\bibitem{CGPS6} 
P. Colli, G. Gilardi, P. Podio-Guidugli, J. Sprekels,
Global existence and uniqueness for a singular/degenerate  
Cahn-Hilliard system with viscosity, 
{\it J. Differential Equations\/} {\bf 254} (2013) 4217-4244.

\bibitem{CGPS7} 
P. Colli, G. Gilardi, P. Podio-Guidugli, J. Sprekels,
Global existence for a strongly coupled 
Cahn-Hilliard system with viscosity, 
{\it Boll. Unione Mat. Ital. (9)} {\bf 5} (2012) 495-513.

{\pier
\bibitem{CGPS8} 
P. Colli, G. Gilardi, P. Podio-Guidugli, J. Sprekels,
Continuous dependence for a nonstandard Cahn-Hilliard 
system with nonlinear atom mobility,
{\it Rend. Semin. Mat. Univ. Politec. Torino\/}
{\bf 70} (2012) 27-52.

\bibitem{CGS1} 
{P. Colli, G. Gilardi, J. Sprekels},
Analysis and optimal boundary 
control of a nonstandard system of phase field equations, 
{\it Milan J. Math.} {\bf 80} (2012) 119-149.

\bibitem{FG} 
E. Fried and M.E. Gurtin, 
Continuum theory of thermally induced phase transitions based on an order 
parameter, {\it Phys. D} {\bf 68} (1993) 326-343.
}

\bibitem{Gurtin} 
M.E. Gurtin, Generalized Ginzburg-Landau and
Cahn-Hilliard equations based on a microforce balance,
{\it Phys.~D\/} {\bf 92} (1996) 178-192.

\bibitem{Podio}
P. Podio-Guidugli, 
Models of phase segregation and diffusion of atomic species on a lattice,
{\it Ric. Mat.} {\bf 55} (2006) 105-118.

\End{thebibliography}

\End{document}

\bye